\documentclass[reqno]{amsproc}
\usepackage{graphic x}
\usepackage{amssymb}
\usepackage[margin=1in]{geometry}
\usepackage{setspace, fullpage}
\usepackage{epstopdf}
\usepackage{amsmath,amsthm,amscd,amssymb}
\usepackage{latexsym}
\usepackage[colorlinks,citecolor=red,pagebackref,hypertexnames=false]{hyperref}
\usepackage{geometry} 
\numberwithin{equation}{section}
\theoremstyle{plain}
\newtheorem{theorem}{{\bf Theorem}}[section]
\newtheorem{lemma}[theorem]{{\bf Lemma}}

\newtheorem{proposition}[theorem]{Proposition}

\theoremstyle{definition}

\newtheorem{case[theorem]}{Case}

\theoremstyle{remark}
\newtheorem{remark}[theorem]{Remark}
\numberwithin{equation}{section}

\begin{document}

\title{\parbox{14cm}{\centering{A Furstenberg-Katznelson-Weiss type theorem on $(d+1)$-point configurations in sets of positive density in finite field geometries}}}

\author{David Covert, Derrick Hart, Alex Iosevich, Steven Senger, and Ignacio Uriarte-Tuero} 

\subitem \email{DavidCovert@mizzou.edu} 
\subitem \email{DNHart@math.rutgers.edu} 
\subitem \email{Iosevich@math.rochester.edu} 
\subitem \email{Senger@math.missouri.edu}
\subitem \email{Ignacio@math.msu.edu}

\thanks{}

\setstretch{1.25}
\begin{abstract} We show that if $E \subset \mathbb{F}_q^d$, the $d$-dimensional vector space over the finite field with $q$ elements, and $|E| \geq \rho q^d$, where $ q^{-\frac{1}{2}}\ll \rho \leq 1$, then $E$ contains an isometric copy of at least $c \rho^{d-1} q^{d+1 \choose 2}$ distinct $(d+1)$-point configurations.
\end{abstract}

\maketitle

\tableofcontents

\section{Background}

 A classical result due to Furstenberg, Katznelson, and Weiss (\cite{FKW90}; see also \cite{B86}) says that if $E \subset {\mathbb R}^2$ has positive upper Lebesgue density, then for any $\delta>0$, the $\delta$-neighborhood of $E$ contains a congruent copy of a sufficiently large dilate of every three point configuration. An example due to Bourgain shows that if the three point configuration in question is an arithmetic progression, then taking a $\delta$-neighborhood is in fact necessary and the result is not otherwise true. However, it seems reasonable to conjecture that if the three point configuration is non-degenerate in the sense that the three points do not lie on the same line, then a set of positive density contains a sufficiently large dilate of this configuration. 

When the size of the point set is smaller than the dimension of the ambient Euclidean space, taking a $\delta$-neighborhood is not necessary, as shown by Bourgain (\cite{B86}). He proves that if $E \subset \mathbb{R}^d$ has positive upper density and $\Delta_k$ is a $k$-simplex (a set of $k+1$ points which spans a $k$-dimensional subspace) with $k<d$, then $E$ contains a rotated and translated image of every large dilate of $\Delta_k$. The cases $k=d$ and $k=d+1$ remain open.  See also, for example, \cite{Berg96}, \cite{F81}, \cite{K07}, \cite{TV06}, and \cite{Z99} on related problems and their connections with discrete analogs. 

In the geometry of the integer lattice $\mathbb{Z}^d$, related problems have been recently investigated by Akos Magyar in \cite{M06} and \cite{M07}. In particular, he proves (\cite{M07}) that if $d>2k+4$ and $E \subset \mathbb{Z}^d$ has positive upper density, then all large (depending on the density of $E$) dilates of a $k$-simplex in $\mathbb{Z}^d$ can be embedded in $E$. Once again, serious difficulties arise when the size of the simplex is sufficiently large with respect to the ambient dimension. 

We aim to investigate an analog of this question in finite field geometries. A step in this direction was taken (\cite{HI07}) by the second and third listed authors. They prove that if $E \subset \mathbb{F}_q^d$, the $d$-dimensional vector space over the finite field with $q$ elements, has $|E| \gtrsim q^{d\frac{k}{k+1}+\frac{k}{2}}$ and $\Delta$ is a $k$-simplex determined by (with vertices lying in) $E$, then there exists $\tau \in \mathbb{F}_q^d$ and $O \in O_d(\mathbb{F}_q)$ such that $\tau+O(\Delta) \subset E$. The result is only non-trivial in the range $d \geq {k + 1 \choose 2}$ as larger simplices are out of range of the methods used.  

Le Anh Vinh has also investigated $k$-point configurations in $\mathbb{F}_q^d$.  He showed in  \cite{Vinh2} that if $E \subset \mathbb{F}_q^d$, $|E| \gtrsim q^{\frac{d-1}{2} + k}$, and $d \geq 2k$  then $E$ contains an isometric copy of every $k$-simplex.  Also, he showed (\cite{Vinh3}) that if an arbitrary set $E \subset \mathbb{F}_q^d$ has size $E \gtrsim q^{\frac{d+2}{2}}$ (for $d \geq 3$), then it determines a positive proportion of all triangles.  Based on an earlier draft of this paper, Vinh proved (\cite{Vinh1}) the $2$-dimensional version of our main theorem (see Theorem \ref{main} below) using graph-theoretic methods.  Namely, if $E \subset \mathbb{F}_q^2$ has size $|E| \geq \rho q^2$ for some $ q^{- \frac{1}{2}} \ll \rho \leq 1$, then the set of triangles determined by $E$ has size $\geq c \rho q^3$.

The purpose of this paper is to address the case of $d$-simplices.  As before we let  $\Delta_k$ denote a $k$-simplex, i.e.~a set of $k+1$ points which span a $k$-dimensional subspace.  Given $E \subset {\mathbb F}_q^d$, let the set of $k$-simplices determined by $E$ up to congruence be denoted by
\[
T_{k}(E)=\{\Delta_k \in E^{k+1}\} \ / \sim
\]
where two $k$-simplices are equivalent if one is a rotated, shifted, reflected copy of the other.

Note that $T_k(E)$ is a natural subset of  $\mathbb{F}_q^{k +1 \choose 2}$ (see Lemma \ref{linalglemma} below).  Our main result is the following. 

\begin{theorem} \label{main}
Let $E \subset \mathbb{F}_q^d$ with $|E| \geq \rho q^d$ for $q^{-1/2} \ll \rho \leq 1$.  Then, there exists $c > 0$ so that 
\[
|T_{d}(E)| \geq c \rho^{d-1} q^{d+1 \choose 2}.
\]
\end{theorem}

\begin{remark}  The viable range for $\rho$ in Theorem \ref{main} is $q^{- (d-\alpha)} \ll \rho \leq 1$, where $\alpha$ is threshold so that
\[
\sum_{x,x^1, \dots , x^d}E(x)E(x^1) \dots E(x^d)S(x-x^1)\dots S(x - x^d) = \frac{|E|^{d+1}}{q^d}(1 + o(1)),
\]
whenever $|E| \gg q^{\alpha}$.  Theorem \ref{dhinges} gives $\alpha = q^{d - \frac{1}{2}}$, although it is reasonable to expect $\alpha = q^\frac{d+1}{2}$
\end{remark}

\begin{remark}
We deal only with finite fields $\mathbb{F}_q$ with characteristic $p > 2$.  We also assume $q$ is much larger than the dimension $d$.  Also, note that the error terms appearing in Theorems \ref{sphere} and \ref{gen} are always of lower order in the effective range of Theorem \ref{main} for $d \geq 2$.
\end{remark}

\begin{remark}
The assumption that $|E| \geq \rho q^d$ implies that the number of $(d+1)$-point configurations determined by $E$ (up to congruence) is at least
\[
\frac{|E|^{d+1}}{\rho q^d \cdot q^{d \choose 2}} \geq \rho^d q^{d+1 \choose 2},
\]
since the size of the subset of the translation group that maps points in $E$ to a set of size $|E|$ is no larger than $ \rho q^d$  and the rotation group is of size $\approx q^{d \choose 2}$.  Our result shaves off a power of $\rho$ from this trivial estimate.

\end{remark}

\section{Proof of the main result (Theorem \ref{main})}

Here, we roughly state the argument.  We prove Theorem \ref{main} by first making a reduction to a statistical statement about hinges (defined below).  Having made this reduction, we next show, using a pigeon-holing argument that for some $x \in E$, the hinge is large.  To finish the argument, we realize a dichotomy.  If the number of transformations mapping the hinge to itself is small, then a purely probabilistic argument gives the number of distinct (incongruent) $(d+1)$-point configurations is what we claim.  If the number of transformations mapping the hinge to itself is large, then a purely combinatorial argument gives the result.  

We start with the statistical reduction.  We observe that if $|E| \geq \rho q^d$, for $\rho$ as above, then it suffices to show that this implies that

\begin{equation}
\left|\left\{ \left(a_{i,j} \right)_{1 \leq i  < j \leq d+1} \in \mathbb{F}_q^{d+1 \choose 2}  : |R_a(E) | > 0\right\}\right| \geq c \rho^{d-1} q^{d +1 \choose 2},
\end{equation}
where 
\[
R_a(E) = \{(y^1, \dots , y^{d+1}) \in E\times \dots \times E : \| y^i - y^j \| = a_{i,j}\},
\]
and
\[
\| x \| = \sum_{j=1}^d x_j^2.
\]

This follows immediately from the following simple linear algebra lemma.  The proof of this lemma will appear in Section \ref{linalgproof} for completeness.

\begin{lemma} \label{linalglemma}
Let $V$ be a simplex with vertices $V_i \in \mathbb{F}_q^d$, where $i = 0, \dots , k$.  Let $W$ be another simplex with vertices $W_i \in \mathbb{F}_q^d$ for $i = 0 , \dots , k$.  Suppose that 

\begin{equation}
\| V_ i  - V_j \| = \| W_i - W_j\|
\end{equation}
for all $i , j$.  Then $ V \sim W$ in the sense of $T_k(E)$.
\end{lemma}

Our main estimate is the following:

\begin{theorem} \label{dhinges}
Suppose that $\alpha_i \in \mathbb{F}_q \backslash \{0\}$ for $i = 1, \dots, d$, and let $E \subset \mathbb{F}_q^d$.  Then,
\[
|\{(x,x^1, \dots , x^d) \in E \times \dots \times E : \| x - x^i  \| = \alpha_i \}| = \frac{|E|^{d+1}}{q^d}(1 +o(1))
\]
whenever $|E| \gg q^{d - \frac{1}{2}}$.
\end{theorem}

This implies that there exists $x \in E$ so that

\begin{equation} \label{d-star}
|\{(x^1, \dots , x^d)\in E \times \dots \times E : \| x - x^i  \| = \alpha_i \}| \geq \frac{|E|^d}{q^d}(1 + o(1)).
\end{equation}

Fix a $d$-tuple $\alpha = (\alpha_i)_{i = 1}^d$, with $\alpha_i \in \mathbb{F}_q \backslash \{ 0 \}$, for $i = 1, \dots , d$.  Define a {\it hinge} $h_{x,\alpha}$ to be the set $\{(x^1, \dots , x^d) \in E \times \dots \times E : \| x - x^i \| = \alpha_i\}$.
Let $M_{x, \alpha} \subset O_d(\mathbb{F}_q)$ denote the set of orthogonal matrices which maps the hinge  $h_{x,\alpha}$ to itself.  We next turn our attention to the following dichotomy:

Suppose that $|M_{x, \alpha}| \leq \rho q^{d \choose 2}$.  By \eqref{d-star}, the number of distinct $d$-point configurations between the $d$ sets $\{x^i \in E : \| x - x^i \| = \alpha_i\}$ is at least
\begin{equation} \label{smallhinge}
\frac{|h_{x, \alpha}|}{|M_{x, \alpha}|}\geq \frac{|E|^d q^{-d} (1 + o(1))}{\rho q^{d \choose 2}}\geq c \rho^{d-1} q^{d \choose 2}.
\end{equation}

We are left only to deal with the case when $|M_{x, \alpha}| > \rho q^{d \choose 2}$.  We put $A_i = \{x^i \in E : \| x - x^i \| = \alpha_i\}.$  It is worthwhile to point out the possibility that $A_i = A_j$. Also, although the sets $A_i$ are not themselves spheres, they are subsets of spheres and therefore inherit some of their intersection properties.  When dealing with the case $|M_{x, \alpha}| > \rho q^{d \choose 2}$ we are faced with two possibilities.  First, suppose that for some $i \in \{ 1, \dots , d \}$ we have that $|A_i| \leq \rho q^{d-1}$.  In this case we utilize the orbit-stabilizer theorem from elementary group theory:

\begin{proposition} \label{ost} $($\cite{Lang}$)$
Let a group $G$ act on a set $S$.  Let $Gs = \{gs : g \in G\}$ be the orbit of $s \in S$, and $G_s = \{g : gs = s\}$ the isotropy group of $s \in S$.  Then there is a bijection between $Gs$ and $G/G_s$.  Consequently,
\[
|Gs| = (G : G_s) = |G| / |G_s|.
\]
\end{proposition}

We let the group $O_d(\mathbb{F}_q)$ act on $\mathbb{F}_q^d$.  Recalling that $|O_d(\mathbb{F}_q)| \approx q ^{d \choose 2}$, and since orthogonal maps preserve the length of a certain vector, we get that the size of the orbit of any point is exactly $q^{d-1}$.  Hence, picking some $z$ from the previously mentioned set $A_i$, we get that the size of the stabilizer group of this element $z$ is
\[
|G_z| = \frac{|G|}{|Gz|} \approx \frac{q^{d \choose 2}}{q^{d-1}}.
\]
The final element here is to notice that
\[
|M_{x, \alpha}| \leq |G_z| |A_i| \leq \frac{q^{d \choose 2}}{q^{d-1}} \cdot \rho q^{d-1} = \rho q^{d \choose 2},
\]
since the number of hinge-preserving orthogonal matrices is no more than the number of orthogonal transformations which fix a given vector $z \in A_i$, times the number of choices for that vector $z$, which is a contradiction.  We may therefore assume $|A_i| > \rho q^{d-1}$ for all $i = 1, \dots , d$.  Recall that we are working with the hinge $h_{x, \alpha} = \{(x^1, \dots , x^d) \in E \times \dots \times E : \| x - x^i \| = \alpha_i\}$, and we aim to show that the number of incongruent $d$-point configurations is bounded below by $c \rho^{d-1} q^{d \choose 2}$.

We start by picking a point $a_1 \in A_1$.  We want to know how many distinct distances occur between $a_1$ and points in the set $A_2$. To achieve this, we count how often a given distance may occur between $a_1$ and the points on $A_2$. This amounts to intersecting $E$ with two spheres: one sphere of a given radius, centered at $a_1$, and the set $A_2$, which is, itself, a sphere intersected with $E$.  The intersection must contain fewer than $q^{d-2}$ possible points on the set $A_2$ which are a given distance from $a_1$. Since $|A_2| > \rho q^{d-1}$, there must be at least $\rho q^{d-1} / q^{d-2} = \rho q$ different distances between $a_1$ and points on $A_2$, by pigeonholing.

For each of the $\rho q$ choices of $a_2$ which are different distances from $a_1$, we need to find the number of 3-point configurations that $a_1$ and $a_2$ can make with points on $A_3$. Now we are intersecting $E$ with spheres of two (possibly the same) radii about $a_1$ and $a_2$ with the sphere containing $S_3$. There can be no more than $q^{d-3}$ points in this intersection, which would each correspond to the same 3-point configuration. So there must be $\rho q^{d-1} / q^{d-3} = \rho q^2$ distinct 3-point configurations for each of the $\rho q$ different pairs we found before, which gives us a total of $\rho q \cdot \rho q^2 = \rho^2 q^3$ different 3-point configurations.  Repeating this process, we see that we will pick up $\rho q^p$ different $(p-1)$-point configurations at each step. If we multiply all of these together, we will get a grand total of 
\begin{equation} \label{largehinge}
\rho q \cdot \rho q^2 \cdot \dots \cdot \rho q^{d-1} =   \rho^{d-1} q^{d \choose 2}
\end{equation} distinct $d$-point configurations.  

From \eqref{smallhinge} and \eqref{largehinge}, we see that in any case, there exist no less than $c \rho^{d-1}q^{d \choose 2}$ many distinct $d$-point configurations.  Since this holds for any fixed vector $\alpha = (\alpha_i)_{i = 1}^d$, and since there are $q - 1$ choices for each $\alpha_i \in \mathbb{F}_q \backslash \{ 0 \}$, then there are at least
\[
c \rho^{d-1} q^{d \choose 2} (q -1)^d \geq c \rho^{d-1} q^{d + 1 \choose 2}
\]
many distinct $(d+1)$-point configurations determined by $E$.

\subsection{Fourier analysis}

\vskip.125in

The Fourier transform of a function $f : \mathbb{F}_q^d \to \mathbb{C}$ is given by
\[
\widehat{f}(m) = q^{-d} \sum_{x \in \mathbb{F}_q^d} f(x) \chi(- x \cdot m)
\]
where $\chi$ is a non-trivial additive character on $\mathbb{F}_q$.  By orthogonality,

\[
 \sum_{x \in \mathbb{F}_q^d} \chi(- x \cdot m) = 
\left\{
\begin{array}{lcc}
q^d  &   &  m = (0,\dots,0) \\
0 &   &   m \neq (0,\dots , 0)
\end{array}
\right.
\]

\begin{lemma} \label{fourierproperties}

Let $f,g : \mathbb{F}_q^d \to \mathbb{C}$.  Then,

\[
\widehat{f}(0,\dots,0) = q^{-d} \sum_{x \in \mathbb{F}_q^d} f(x),
\]

\[
q^{-d} \sum_{x \in \mathbb{F}_q^d} f(x) \overline{g(x)} = \sum_{m \in \mathbb{F}_q^d} \widehat{f}(m) \overline{\widehat{g}(m)},
\]

\[
f(x) = \sum_{m \in \mathbb{F}_q^d} \widehat{f}(m) \chi(x \cdot m).
\]

\end{lemma}

\section{Proof of Theorem \ref{dhinges}} 
In order to prove Theorem \ref{dhinges} we will actually prove the more general following theorem. \begin{theorem} \label{gen}
Let $r>2$ be an integer, and let  $H_{r,\alpha}$ represent the set of $r-$hinges, with distances $\alpha= \{\alpha_i\}_{i=1}^{r-1}$, which are present in $E$.  That is,
\[
H_{r,\alpha} = \{(x,x^1, \dots x^{r-1}) \in E \times \dots \times E : \| x - x^i \| = \alpha_i\},
\]
where $\alpha_i \neq 0$ for $i = 1 , \dots , r-1$.  Then, 
\[
|H_{r,\alpha}| = \frac{|E|^{r}}{q^{r-1}}(1 + o(1)),\] whenever
$|E|\gg q^{\frac{2r-5}{2r-4}d+\frac{1}{2r-4}}$
\end{theorem}

Setting $r=d+1$ in Theorem \ref{gen} gives Theorem \ref{dhinges}.

We will need the following estimates whose proof we delay to the end of the paper.
\begin{theorem} \label{sphere}
Let $S_t = \{x \in \mathbb{F}_q^d : \| x \| = t\}$.  Identify $S_t$ with its characteristic function.  For $t \neq 0$,
\begin{equation} \label{sizeofsphere}
|S_t| = q^{d-1}(1 + o(1))
\end{equation}
and if also $m \neq (0, \dots , 0)$,
\begin{equation} \label{decaysphere}
| \widehat{S}_t(m) | \lesssim q^{- \frac{d+1}{2}}.
\end{equation}
\end{theorem}

We will proceed by induction.  Before we handle the case $r=3$ we first observe the following estimate which originally appeared in \cite{IR}.  

\begin{lemma} \label{2hinge}
Using the notation as above, we have $|H_{2,\alpha}| = \frac{|E|^2}{q} + O(q^{\frac{d-1}{2}}|E|)$.
\end{lemma}

To see this, write
\begin{align*}
|H_{2,\alpha}| &= \sum_{x,y} E(x) E(y) S(x - y)
\\
&= q^{2d} \sum_m \left| \widehat{E}(m) \right|^2 \widehat{S}(m)
\\
&= q^{-d} |E|^2 |S| + q^{2d} \sum_{m \neq 0} \left| \widehat{E}(m) \right|^2 \widehat{S}(m)
\end{align*}
and
\begin{align*}
q^{2d} \left| \sum_{m \neq 0} \left| \widehat{E}(m) \right|^2 \widehat{S}(m)\right| &\leq 2 q^{2d} q^{- \frac{d+1}{2}} q^{-d} |E|
\\
&= 2q^{\frac{d-1}{2}}|E|.
\end{align*}

We now illustrate the base step.  First we write
\[
|H_{3,\alpha}|=\sum_{x\in E} |E \cap (x - S)| ^2.
\]
Now,
\[
|E \cap (x - S)|=\sum_y E(y)S(x-y)=q^d\sum_m \widehat{E}(m)\widehat{S}(m)\chi(m\cdot x)
\]
\[
=|E||S|q^{-d}+q^d\sum_{m \neq 0} \widehat{E}(m)\widehat{S}(m)\chi(m\cdot x),
\]
which gives

\begin{align*}
|H_{3,\alpha}| &= \sum_{x\in E} |E \cap (x - S)| ^2 
\\
&= |E|^3|S|^2q^{-2d}+2|E||S|q^d\sum_{m\neq 0} |\widehat{E}(m)|^2|\widehat{S}(m)|
+q^{2d}\sum_{x}\left|\sum_{m\neq 0} \widehat{E}(m)\widehat{S}(m)\chi(m\cdot x)\right|^2
\\
&=|E|^3|S|^2q^{-2d}+O\left(|E|^2|S|q^{-d}q^{(d-1)/2}+q^{3d}\sum_{m \neq 0} |\widehat{E}(m)|^2|\widehat{S}(m)|^2\right)
\\
&= |E|^3|S|^2q^{-2d}+O\left(|E|^2|S|q^{-d}q^{(d-1)/2}+q^{d-1}|E|\right).
\end{align*}

If $|E|\gg q^{\frac{d+1}{2}}$ then 
\[
|H_{3,\alpha}|= |E|^3 q^{-2}(1+o(1)),
\]

For the inductive step, assume that we are in the case $|H_{r,\alpha}| = \frac{|E|^r}{q^{r-1}}(1 + o(1))$ for $|E|\gg q^{\frac{2r-5}{2r-4}d+\frac{1}{2r-4}}.$
We begin by writing
\begin{align*}
|H_{r+1,\alpha}| &= \sum_{x, x^1 , \dots, x^r} H_{r,\alpha}(x , x^1 , \dots , x^{r-1}) E(x^r) S(x - x^r)
\\
&= q^{(r+1)d} \sum_{m} \widehat{H}_{r,\alpha}(m, 0 , \dots , 0) \widehat{S}(m) \widehat{E}(m)
\\
&= q^{-d} |E| |S| |H_{r,\alpha}| + q^{(r+1)d} \sum_{m \neq 0}  \widehat{H}_{r,\alpha}(m, 0 , \dots , 0) \widehat{S}(m) \widehat{E}(m)
\\
&= q^{-d} |E| |S| |H_{r,\alpha}| + R.
\end{align*}

Applying Cauchy-Schwarz gives

\begin{align*}
R^2 &\leq q^{2d(r+1)} \sum_{m \neq 0} |\widehat{S}(m)|^2 |\widehat{E}(m)|^2 \sum_{m \neq 0} |\widehat{H}_{r,\alpha}(m, 0 , \dots, 0)|^2
\\
&\lesssim q^{2d(r+1)}q^{-d-1}q^{-d}|E| \sum_{m \neq 0} |\widehat{H}_{r,\alpha}(m, 0 , \dots , 0)|^2
\\
&\leq q^{2dr - 1}|E| \sum_{m} |\widehat{H}_{r,\alpha}(m, 0, \dots , 0)|^2
\end{align*}

Also,  we have that 

\begin{align*}
\widehat{H}_{r,\alpha}(m, 0 , \dots , 0) &= q^{-rd} \sum_{x, x^1, \dots , x^{r-1}} \chi(x \cdot m) E(x) E(x^1) \dots E(x^{r-1}) S(x - x^1) \dots S(x - x^{r-1})
\\
&= q^{-rd + d} \widehat{f}(m)
\end{align*}

where 
\[
f(x) = E(x) \sum E(x^1) \dots , E(x^{r-1}) S(x - x^1) \dots S(x - x^{r-1}) = E(x)|E \cap (x - S)|^{r-1}.
\]

Since $|E \cap (x - S)| \leq q^{d-1}$, it follows that

\begin{align*}
A = \sum_{m} |\widehat{H}_{r,\alpha}(m, 0 , \dots , 0 ) |^2 &= q^{-2rd + 2d} \sum_{m} |\widehat{f}(m)|^2
\\
&= q^{-2rd + d} \sum_{x} |f(x)|^2
\\
&\leq q^{-2rd + d} \left( q^{d-1}\right)^{2(r-2)} |H_{3,\alpha}|
\end{align*}
and 
\[
A \lesssim q^{-2rd + d} \left( q^{d-1}\right)^{2(r-2)}  |E|^3 q^{-2}(1+o(1)).
\]
Finally,
\[
R^2 \lesssim q^{ - 3}q^{ d} \left( q^{d-1}\right)^{2(r-2)}  |E|^4 (1+o(1))\leq  q^{(2r-3)d-2r+1}|E|^4 (1+o(1)).
\]

Therefore, 
\[
|H_{r+1,\alpha}| = q^{-d} |E| |S| |H_{r,\alpha}| + O(q^{d \frac{2r-3}{2} - r + \frac{1}{2}}|E|^2),
\]
and we have that 
\[
|H_{r+1,\alpha}| = \frac{|E|^{r+1}}{q^r}(1 + o(1)),
\]
whenever
\[
|E|\gg q^{\frac{2r-3}{2r-2}d+\frac{1}{2r-2}}.
\]

\section{Proof of Lemma \ref{linalglemma}} \label{linalgproof}

Let $\pi_r(x)$ denote the $r$-th coordinate of $x$.  By translating, we may assume that $V_0 = \vec{0}$.  We may also assume that $V_1 , \dots , V_k$ are contained in $\mathbb{F}_q^k$.  The condition that $\| V_i - V_j \| = \| W_i - W_j\|$ for all $ i, j$ implies that

\begin{equation} \label{fixeddistance}
\sum_{r = 1}^k \pi_r(V_i)\pi_r(V_j) = \sum_{r=1}^k \pi_r(W_i) \pi_r (W_j).
\end{equation}

Let $T$ be the transformation uniquely defined by $T(V_i) = W_i$.  To show that $T$ is orthogonal it suffices to show that $\| Tx \| = \| x \|$ for all $x$.  By assumption, the $V_i$'s form a basis, so we have
\[
x = \sum_{i} t_i V_i.
\]
Thus, by \eqref{fixeddistance}, we have that
\[
\| Tx \| = \sum_{r} \sum_{i , j} t_i t_j \pi_r(W_i) \pi_r(W_j) = \sum_r \sum_{i,j} t_i t_j \pi_r(V_i) \pi_r(V_j) = \| x \|,
\]
giving the result.

\section{Proof of Theorem \ref{sphere}} \label{sphereproof}

For any $l\in{\mathbb F}^d_q$, we have
\begin{equation} \label{sphereparade}
\begin{array}{llllll} \widehat{S}_t(l)&=&
q^{-d} \displaystyle\sum_{x \in {\mathbb F}^d_q} q^{-1} \sum_{j \in {\mathbb F}_q} \chi(
j(\|x\|-t)) \chi( -  x \cdot l)\\ \hfill \\&=&q^{-1}\delta(l) +
q^{-d-1} \displaystyle\sum_{j \in {\mathbb F}^{*}_q} \chi(-jt) \sum_{x}
\chi( j\|x\|) \chi(- x \cdot l),\\
\end{array}\end{equation}where the notation $\delta(l)=1$ if $l=(0\ldots,0)$ and $\delta(l)=0$ otherwise. 

Now
\[
\widehat{S}_t(l)=q^{-1}\delta(l)+   Q^d q^{-\frac{d+2}{2}} \sum_{j \in {\mathbb F}^{*}_q}
\chi\left(\frac{\|l\|}{4j}+jt\right)\eta^d(-j).
\]
In the last line we have completed the square, changed $j$ to
$-j$, and used $d$ times the Gauss sum equality
\begin{equation} \sum_{c\in {\mathbb F}_q} \chi(jc^2) = \eta(j)\sum_{c\in{\mathbb F}_q}\eta(c)\chi(c)=\eta(j)\sum_{c\in{\mathbb F}_q^*}\eta(c)\chi(c) =Q\sqrt{q}\,\eta(j), \label{gauss}\end{equation} where the constant $Q$ equals $\pm1$ or $\pm i$, depending on $q$, and $\eta$ is the quadratic multiplicative character (or the Legendre symbol) of ${\mathbb F}_q^*$.  (see, e.g. \cite{LN97}, for more information).

The conclusion to both parts of Theorem \ref{sphere} now follows from the following classical estimate due to A. Weil (\cite{Weil}).
\begin{theorem} \label{kloosterman} Let
\[
K(a)=\sum_{s \not=0} \chi(as+s^{-1}) \psi(s),
\]
 where $\psi$ is a multiplicative character on ${\mathbb F}_q\backslash\{0\}$. Then if $a \neq 0$,
\[
|K(a)| \leq 2 \sqrt{q}.
\]
\end{theorem}

\end{document}